\documentclass[11pt,a4paper]{article}
\usepackage[T1]{fontenc}
\usepackage{lmodern}
\usepackage[margin=28mm]{geometry}
\usepackage{amsmath,amssymb,amsthm,mathtools}
\usepackage{microtype,enumitem,needspace}
\usepackage[dvipsnames]{xcolor}
\usepackage[colorlinks=true,linkcolor=MidnightBlue,citecolor=MidnightBlue,urlcolor=MidnightBlue]{hyperref}
\hypersetup{pdftitle={Axioms for relative entropy},pdfauthor={Roberto Longo}}
\setlist[enumerate]{leftmargin=*,itemsep=5pt,topsep=5pt}
\newtheorem{theorem}{Theorem}[section]
\newtheorem{lemma}[theorem]{Lemma}
\newtheorem{proposition}[theorem]{Proposition}
\theoremstyle{definition}
\newtheorem{remark}[theorem]{Remark}
\newcommand{\Sn}{\mathcal S_{\mathrm n}}
\newcommand{\vN}{von Neumann}
\newcommand{\ot}{\,\otimes\,}
\newcommand{\relsep}{\mathrel{|\!|}}
\newcommand{\Dmax}[1][]{D^{\max#1}}
\newcommand{\id}{\operatorname{id}}
\newcommand{\supp}{\operatorname{s}}
\newcommand{\tr}{\operatorname{Tr}}
\newcommand{\htr}{\operatorname{tr}}
\newcommand{\C}{\mathbb C}

\newcommand{\step}[1]{\par\smallskip\noindent\textit{#1}\quad}
\title{\huge Axioms for relative entropy}
\author{Roberto Longo\\[3pt]\small Department of Mathematics, University of Rome Tor Vergata\\\small Via della Ricerca Scientifica, 1 - 00133 Rome, Italy}
\date{Revised 14 November 2026}

\begin{document}
\maketitle 
\begin{abstract}
For normal states of arbitrary \vN\ algebras, we characterize
positive finite multiples of Araki relative entropy among
non-degenerate divergences. The axioms are monotonicity under normal
unital completely positive maps, tensor additivity, superadditivity
for product second states, and separate lower semicontinuity in the
predual norm.
Alternatively, a
conditional-expectation identity replaces the two tensor assumptions.
The proof uses order bounds and asymptotic estimates on tensor powers.
An appendix proves the known analytic estimates needed here, using
standard modular theory and Haagerup $L^p$ spaces. We also obtain the
corresponding characterization 
in the C*-algebraic setting. 
\end{abstract}

\section*{Introduction}
Relative entropy is a primary information-theoretic notion. It measures the distinguishability of two states of a system. 
It was introduced by Kullback and Leibler \cite{KL} for classical probability distributions. If $p = (p_1, \dots, p_r)$, $q= (q_1, \dots , q_r)$ are finite probability measures, their relative entropy is given by
\[
 S_{\C^r}(p\relsep q)= \sum_i p_i\log\frac{p_i}{q_i}\, .
\]
Umegaki introduced the relative entropy in the semifinite setting
\cite{Umegaki62}; for density matrices of a matrix algebra $A$,
\[
S_A(\rho\relsep\sigma)=\tr\rho(\log\rho-\log\sigma).
\]
Araki extended it to normal positive functionals
on arbitrary \vN\ algebras \cite{Araki76,Araki77}, thus for arbitrary quantum systems. For normal faithful states
$\varphi,\psi$ of a \vN\ algebra $M$ in standard form, Araki's
definition is
\[
 S_M(\varphi\relsep\psi)
 =-\langle\xi_\varphi,
          \log\Delta_{\psi,\varphi}\,\xi_\varphi\rangle.
\]
Here $\Delta_{\psi,\varphi}$ is the relative
modular operator between cyclic and separating vectors $\xi_\psi$, $\xi_\varphi$ representing $\psi$, $\varphi$ .
The expression is interpreted
as an extended spectral integral. 

Uhlmann's construction by interpolation of quadratic forms
\cite{Uhlmann77} agrees with Araki's for normal functionals; see
Hiai--Ohya--Tsukada \cite[Section~3]{HOT83}. 
Kosaki's 
formula \cite{Kosaki} gives another important variational description of the same functional.
A recent, related description as $f$-divergence is given by the integral representation of
Correa da Silva--Fr\"ob--Lechner--Sangaletti
\cite[Theorem~3.19]{CFLS2026}. 

The author's interest in relative entropy originated from work on
von Neumann algebras and quantum black holes \cite{Longo97}. Quantum
field theory provides a natural setting for infinite quantum systems,
in which relative entropy plays a crucial role, in particular in the
study of local algebras, entropy-energy inequalities and black hole
entropy. See \cite{HollandsLongo25QNEC,Longo25Bekenstein,LongoXu18}
for examples of recent developments and further references.

The question here is whether relative entropy is determined by a small
number of structural properties. Petz \cite{Petz92} and Ohya--Petz
\cite{OP} gave earlier matrix-algebra characterizations involving
conditional expectations. Donald \cite{Donald86} gives a distinct axiomatic construction of a
generalized relative entropy and proves its agreement with the standard
relative entropy for normal state pairs on injective \vN\ algebras.
For nuclear C*-algebras, Petz
\cite[Theorem~1.2]{Petz1999} states a characterization using the
conditional-expectation identity, monotonicity under unital Schwarz
maps, a direct-sum rule, and lower semicontinuity.

The asymptotic argument in this paper also has antecedents. Matsumoto
\cite{Matsumoto} characterizes finite-dimensional relative entropy,
up to scale, by monotonicity, additivity on tensor powers, and
asymptotic lower semicontinuity. Wilming--Gallego--Eisert \cite{WGE}
obtain a characterization for matrix states with a faithful second
state from continuity in the first state, monotonicity, additivity,
and superadditivity. Their use of superadditivity to pass from
regularity on one algebra to an estimate on its tensor powers is the
mechanism used in Lemma~\ref{lem:alsc}, with lower semicontinuity in
place of continuity.

In this paper we give a characterization for arbitrary \vN\ algebras,
without injectivity or faithfulness assumptions. Starting with a
divergence $D$, that is a family of maps $D_M$ taking values in
$[0,+\infty]$ and defined for all normal state pairs on arbitrary
\vN\ algebras $M$, we show that $D$ is a positive finite multiple of
Araki relative entropy, with the same constant for every algebra and
every pair of states, if $D$ satisfies four axioms: monotonicity,
tensor additivity, superadditivity for product second states, and
separate norm lower semicontinuity.
 We assume that the
divergence $D$ is non-degenerate, meaning that it has a positive finite value at
one state pair on one algebra. Neither the location of that pair nor
its value is prescribed. 
The formula on finite commutative algebras is
derived. 
Alternatively, the conditional-expectation identity can replace
tensor additivity and superadditivity, with monotonicity and separate
norm lower semicontinuity retained.
An analogous characterization holds in the C*-algebraic setting. 
The C*-algebra proof uses the minimal tensor product and Kaplansky
density, without a nuclearity assumption. 

The analytic estimates in Proposition~\ref{prop:aep} are known:
they follow from Fawzi--Gao--Rahaman \cite{FGR}, with the estimates
for separation of tensor powers by positive contractions due to
Pautrat--Wang \cite{PW}. Appendix~\ref{app:aep} gives a pedagogical
proof of the dominated case needed here, using standard modular
theory and Haagerup $L^p$ spaces \cite{HaagerupLp}. Its role is
expository; it supplies established analytic input independently of
the function $D$ being characterized.

Sections~1--4 state the axioms and prove the \vN\ result.
Section~\ref{sec:alternative} gives the conditional-expectation
formulation, and Section~\ref{sec:cstar} proves the C*-algebra result.

\section{Setting, axioms, and plan of the proof}

\subsection{States, norms, and unital completely positive maps}
Let $M$ be a \vN\ algebra, $M_*$ its predual, and $\Sn(M)$ its normal
state space. Tensor products of \vN\ algebras are denoted by $M\otimes N$.
As usual, the support $\supp(\varphi)$ is the least projection $e\in M$ with
$\varphi(e)=1$. Faithfulness means $\supp(\varphi)=1$. We call $M$
$\sigma$-finite if it has a faithful normal state.

We consider a family of maps $D = (D_M)_M$
\[
 D_M:\Sn(M)\times\Sn(M)\longrightarrow[0,+\infty],
\]
indexed by all von Neumann algebras $M$. 
No symmetry, bilinearity, or continuity is included in this
definition. The value $+\infty$ is allowed and is important below.
We call the family $D$ a \emph{divergence}, and we say that the divergence $D$ is non-degenerate if there exist a
\vN\ algebra $M_0$ and normal states $\varphi_0,\psi_0$ such that
\begin{equation}\label{eq:proper}
 0<D_{M_0}(\varphi_0\relsep\psi_0)<+\infty.
\end{equation}
Non-degeneracy is a condition on the family: such a pair need not exist
on each algebra. In particular, no positive value is required on
$\C$, which has only one state. Vanishing on equal states is not
part of the definition; it will follow from non-degeneracy and the axioms.

A normal unital completely positive map (UCP) $\Phi:N\to M$ induces the
map $\varphi\mapsto\varphi\circ\Phi$ from $\Sn(M)$ to $\Sn(N)$.
For finite commutative algebras we often describe this induced map
directly on probability vectors. If $T_{ji}\geq0$ and
$\sum_jT_{ji}=1$, its action is $p\mapsto Tp$, and its adjoint is
\[
 \Phi_T:\C^J\to\C^I,\qquad
 (\Phi_T f)_i=\sum_jT_{ji}f_j.
\]
This adjoint is unital and completely positive. Consequently, the monotonicity axiom applies to every stochastic map used in the proof.

\subsection{The entropy to be characterized}
We use the nonnegative convention for Araki relative entropy
$S_M(\varphi\relsep\psi)$ \cite{Araki76,Araki77,OP}. In standard form, let
$\xi_\varphi$ be the natural-cone vector representing $\varphi$, and let
$\Delta_{\psi,\varphi}$ be the relative modular operator in this order.
Then
\begin{equation}\label{eq:araki}
 S_M(\varphi\relsep\psi)
 =-\langle\xi_\varphi,\log\Delta_{\psi,\varphi}\,\xi_\varphi\rangle.
\end{equation}
The expression is interpreted as an extended spectral integral and is well-defined. 

\begin{theorem}[Axiomatic characterization]\label{thm:main}
Let $D=(D_M)_M$ be a non-degenerate divergence on all \vN\ algebras and all
their normal state pairs. Assume the following four conditions.
\begin{enumerate}[label=\textup{(A\arabic*)},ref=A\arabic*]
\item\label{ax:dpi}\textbf{Monotonicity.}
For every normal unital completely positive map $\Phi:N\to M$,
\[
 D_N(\varphi\circ\Phi\relsep\psi\circ\Phi)
 \leq D_M(\varphi\relsep\psi).
\]

\item\label{ax:add}\textbf{Tensor additivity.}
For $\varphi_i,\psi_i\in\Sn(M_i)$,
\[
 D_{M_1\ot M_2}(\varphi_1\otimes\varphi_2\relsep\psi_1\otimes\psi_2)
 =D_{M_1}(\varphi_1\relsep\psi_1)+D_{M_2}(\varphi_2\relsep\psi_2).
\]

\item\label{ax:super}\textbf{Superadditivity for product second states.}
For $\omega\in\Sn(M\ot N)$, denote its marginal states by
$\omega_M(x)=\omega(x\otimes1)$ and
$\omega_N(y)=\omega(1\otimes y)$. Then
\[
 D_{M\ot N}(\omega\relsep\psi_M\otimes\psi_N)
 \geq D_M(\omega_M\relsep\psi_M)+D_N(\omega_N\relsep\psi_N).
\]

\item\label{ax:lsc}\textbf{Separate lower semicontinuity in norm.}
With either state fixed, the function is lower semicontinuous in the
other state for the predual norm. Explicitly,
\begin{align*}
 \|\varphi_k-\varphi\|\to0
 &\quad\Longrightarrow\quad
 D_M(\varphi\relsep\psi)
 \leq\liminf_kD_M(\varphi_k\relsep\psi),\\
 \|\psi_k-\psi\|\to0
 &\quad\Longrightarrow\quad
 D_M(\varphi\relsep\psi)
 \leq\liminf_kD_M(\varphi\relsep\psi_k).
\end{align*}
\end{enumerate}
Then there is a constant $0<\lambda<\infty$, independent of $M$
and of the states, such that
\[
 D_M(\varphi\relsep\psi)=\lambda S_M(\varphi\relsep\psi)
\]
for every normal state pair on every \vN\ algebra. More precisely,
for $\delta=(1,0)$ and $u=(\tfrac12,\tfrac12)$,
\begin{equation}\label{eq:scale}
 0<b:=D_{\C^2}(\delta\relsep u)<\infty,
 \qquad \lambda=\frac{b}{\log2}.
\end{equation}
Conversely, every positive finite multiple of Araki relative entropy
 satisfies these four conditions.
\end{theorem}

Note that normal
UCP maps between different algebras are permitted. The theorem has no
injectivity, hyperfiniteness, or faithfulness hypothesis. It implies the
formula on finite commutative algebras
\[
 D_{\C^r}(p\relsep q)=\lambda\sum_i p_i\log\frac{p_i}{q_i},
\]
with $0\log(0/q)=0$ and $p\log(p/0)=+\infty$ for $p>0$.

\begin{remark}[What fixes the scale]\label{rem:scale}
Non-degeneracy does not fix a numerical scale. Lemma~\ref{lem:proper}
derives the positive finite value $b$ in \eqref{eq:scale} from the
four axioms. The optional condition
\begin{equation}\label{eq:calibration}
 D_{\C^2}(\delta\relsep u)=\log2
\end{equation}
then selects $D=S$.
The degenerate divergences are the identically zero assignment, the identically
infinite assignment, and the assignment that is zero on equal states
and infinite on distinct states; all satisfy the four axioms. 
\end{remark}

\subsection{How the proof is organized}
First we show, without using any entropy characterization, that
non-degeneracy forces $0<b<\infty$. We then divide $D$ by $b/\log2$.
For this normalized divergence the proof has two comparisons.
Order bounds obtained from UCP maps
into $\C^2$, together with approximation on tensor powers, give
$D\leq S$. Positive contractions separating the tensor powers,
together with superadditivity, give $D\geq S$. These comparisons are first made for
pairs with an order bound $\varphi\leq C\psi$, $C<\infty$.
Replacing $\psi$ by $(1-t)\psi+t\varphi$ produces such a bound for
every $t>0$; monotonicity and lower semicontinuity then let $t$ tend
to zero. Finally a support-corner argument removes $\sigma$-finiteness.

Section~2 derives the binary value from non-degeneracy, introduces this
normalization, and proves the elementary consequences used in both
comparisons. Section~3 states the asymptotic facts independently of the
unknown function $D$. Section~4 proves uniqueness.
Section~\ref{sec:alternative} gives the conditional-expectation
formulation, and Section~\ref{sec:cstar} proves the C*-algebra extension.
Appendix~\ref{app:aep} reproves the known asymptotic statements, using
standard modular and noncommutative $L^p$ background. Its role is to
make the analytic input self-contained at that level; it does not use
the characterization theorem that it helps prove.

\section{Non-degeneracy, the scale, and elementary consequences}

An elementary observation will be used repeatedly. If precomposition
with normal UCP maps carries two state pairs to each other, monotonicity in
both directions makes their values under $D$ equal. The compositions
need only fix the two specified states; they need not fix all states.

\begin{lemma}[Non-degeneracy determines a positive finite binary value]
\label{lem:proper}
Let $D$ be non-degenerate and satisfy \ref{ax:dpi}--\ref{ax:lsc}. Then
$D_M(\chi\relsep\chi)=0$ for every normal state $\chi$, and
\[
 0<b:=D_{\C^2}(\delta\relsep u)<\infty.
\]
\end{lemma}
\begin{proof}
\step{Step 1: all diagonal values vanish.}
Let $\omega$ be the unique state on $\C$. Scalar inclusion and state
evaluation are normal UCP maps in opposite directions. Monotonicity
therefore gives, for every normal state $\chi$,
\[
 D_M(\chi\relsep\chi)=D_\C(\omega\relsep\omega)=:a.
\]
Scalar inclusion also shows $a\leq D_M(\varphi\relsep\psi)$ for every
pair. The positive finite value in \eqref{eq:proper} implies that $a$
is finite. Tensor additivity gives $a=2a$, so $a=0$. In particular,
the pair witnessing non-degeneracy consists of distinct states.

\step{Step 2: the binary value cannot be zero.}
Suppose $b=0$. Tensor additivity gives zero on
$(\delta^{\otimes k},u^{\otimes k})$. Recording whether the outcome
is $0^k$ is a stochastic map, so monotonicity implies
\[
 D_{\C^2}\bigl(\delta\relsep(2^{-k},1-2^{-k})\bigr)=0,
 \qquad k\geq1.
\]
If $\varphi\leq C\psi$ for a finite $C$, choose $k\geq1$ with
$2^k\geq C$ and set
\[
 \tau=\frac{\psi-2^{-k}\varphi}{1-2^{-k}}.
\]
This is a normal state. Precomposition with the normal UCP map
$x\mapsto(\varphi(x),\tau(x))$ sends
$(\delta,(2^{-k},1-2^{-k}))$ to $(\varphi,\psi)$. Hence
$D_M(\varphi\relsep\psi)=0$ for every dominated pair.

For arbitrary states put $\psi_t=(1-t)\psi+t\varphi$, $0<t<1$.
Since $\varphi\leq t^{-1}\psi_t$, the preceding conclusion gives
$D_M(\varphi\relsep\psi_t)=0$. As $\psi_t\to\psi$ in norm, lower
semicontinuity in the second argument yields
$D_M(\varphi\relsep\psi)=0$. This contradicts non-degeneracy. Thus $b>0$.

\step{Step 3: reduce a finite value to distinct binary distributions.}
Suppose now that $b=+\infty$, and choose the pair
$(\varphi_0,\psi_0)$ in \eqref{eq:proper}. Since these states are
distinct, rescaling a self-adjoint element that separates them gives
a positive contraction $e\in M_0$ with
\[
 \tfrac14 1\leq e\leq\tfrac34 1,
 \qquad \varphi_0(e)\ne\psi_0(e).
\]
The normal UCP map $(z_0,z_1)\mapsto z_0e+z_1(1-e)$ gives distinct,
strictly positive binary distributions
\[
 p=(\varphi_0(e),1-\varphi_0(e)),\qquad
 q=(\psi_0(e),1-\psi_0(e)).
\]
Monotonicity implies
\[
 d:=D_{\C^2}(p\relsep q)
 \leq D_{M_0}(\varphi_0\relsep\psi_0)<\infty.
\]
We will show that this is incompatible with $b=+\infty$.

\step{Step 4: construct many binary marginals close to $\delta$.}
Define the ordinary number
\[
 s=\sum_{i=0}^1p_i\log\frac{p_i}{q_i}>0,
\]
whose strict positivity follows from $p\ne q$. 
Fix $0<R<s$ and put
$P_n=p^{\otimes n}$, $Q_n=q^{\otimes n}$. On $\{0,1\}^n$ let
\[
 A_n=\left\{x:\log\frac{P_n(x)}{Q_n(x)}\geq nR\right\},
 \quad \alpha_n=P_n(A_n),\quad\beta_n=Q_n(A_n).
\]
Under $P_n$, the logarithm is a sum of independent identically
distributed random variables with mean $s$ and finite variance.
The weak law of large numbers gives $\alpha_n\to1$. On $A_n$ one
has $Q_n(x)\leq e^{-nR}P_n(x)$, so $\beta_n\leq e^{-nR}$.

Set
\[
 m_n=\left\lfloor\frac{nR}{\log2}\right\rfloor,
 \quad a_n=\delta^{\otimes m_n},\quad v_n=u^{\otimes m_n}.
\]
For large $n$, $m_n\geq1$ and $\beta_n\leq2^{-m_n}$. Thus
\[
 r_n=\frac{v_n-\beta_na_n}{1-\beta_n}
\]
is a probability distribution. The stochastic map which outputs
$a_n$ on $A_n$ and $r_n$ on its complement sends $Q_n$ to $v_n$
and $P_n$ to
\[
 \omega_n=\alpha_na_n+(1-\alpha_n)r_n.
\]
Consequently, monotonicity and tensor additivity give
\[
 nd\geq D_{\C^{2^{m_n}}}(\omega_n\relsep v_n).
\]
Moreover $\|\omega_n-a_n\|_1\leq2(1-\alpha_n)\to0$, so every
binary marginal $\omega_{n,j}$ tends to $\delta$, uniformly in $j$.

\Needspace{8\baselineskip}
\step{Step 5: use the infinite binary value to contradict finiteness.}
Fix any finite $L>0$. Since $D_{\C^2}(\delta\relsep u)=+\infty$,
lower semicontinuity in the first argument gives a norm neighborhood
of $\delta$ on which $D_{\C^2}(\chi\relsep u)>L$. For all large $n$
every $\omega_{n,j}$ belongs to this neighborhood. Superadditivity
then yields
\[
 nd\geq\sum_{j=1}^{m_n}D_{\C^2}(\omega_{n,j}\relsep u)
 \geq m_nL.
\]
Dividing by $n$ and using $m_n/n\to R/\log2$ gives
$d\geq (R/\log2)L$. Since $L$ is arbitrary, this contradicts
$d<\infty$. Hence $b<\infty$, completing the proof.
\end{proof}

The lemma uses only the four axioms and elementary probability on
finite sets. In particular, neither Araki's entropy nor the asymptotic
results of Section~3 enter its proof.
Define
\begin{equation}\label{eq:rescaling}
 \widehat D_M=\frac{\log2}{b}D_M.
\end{equation}
This is again non-degenerate and satisfies the same four axioms, with
$\widehat D_{\C^2}(\delta\relsep u)=\log2$.
For the remainder of this section and in Section~4 we denote
$\widehat D$ simply by $D$. We will prove $\widehat D=S$ and then
undo \eqref{eq:rescaling}. Thus the following normalization is a
reduction of the proof, not an additional axiom.

\begin{lemma}[Values on $\C^2$]\label{lem:binary}
For the normalized divergence just introduced, \ref{ax:dpi} and
\ref{ax:add} imply
$D_M(\chi\relsep\chi)=0$ for every normal state $\chi$, and
\[
 D_{\C^2}\bigl(\delta\relsep(t,1-t)\bigr)=-\log t,
 \qquad 0<t\leq1.
\]
The value at $t=0$ is $+\infty$.
\end{lemma}
\begin{proof}
\step{Step 1: identify and evaluate all diagonal values.}
Let $\omega$ be the unique state on $\C$. For a normal state $\chi$ of
$M$, the maps
\[
 \iota:\C\to M,\quad\iota(z)=z1,
 \qquad E_\chi:M\to\C,\quad E_\chi(x)=\chi(x)
\]
are normal unital completely positive. Since
$\chi\circ\iota=\omega$ and $\omega\circ E_\chi=\chi$,
\ref{ax:dpi} gives
\[
 D_\C(\omega\relsep\omega)
 \leq D_M(\chi\relsep\chi)
 \leq D_\C(\omega\relsep\omega).
\]
Thus every diagonal value is the same number, say $a$.
Discarding the outcome of the normalized pair $(\delta,u)$ gives
$a\leq\log2$, so $a$ is finite.
The diagonal value at $(\chi\otimes\chi,\chi\otimes\chi)$ is also $a$,
whereas \ref{ax:add} says it is $2a$. Hence $a=0$.
The finiteness just established is needed here: the extended equality
$+\infty=2(+\infty)$ would give no conclusion.

\step{Step 2: show monotonicity of the binary function.}
Write $q_t=(t,1-t)$ and $g(t)=D_{\C^2}(\delta\relsep q_t)$ for
$0\leq t\leq1$. If $0\leq s\leq t\leq1$ and $s<1$, put
\[
 r=\frac{t-s}{1-s},\qquad
 T_{s,t}=\begin{pmatrix}1&r\\0&1-r\end{pmatrix}.
\]
The entries are nonnegative and each column sums to one, so this is a
stochastic map. Direct multiplication gives
\[
 T_{s,t}\delta=\delta,
 \qquad T_{s,t}q_s=(s+r(1-s),(1-r)(1-s))=q_t.
\]
Monotonicity therefore implies $g(t)\leq g(s)$.
Thus $g$ is nonincreasing. Step~1 also gives $g(1)=0$.

\step{Step 3: compare the tensor-product pair with a pair on $\C^2$.}
Fix $0<t<1$ and $n\geq1$. On $\Omega_n=\{0,1\}^n$, let
\[
 P_n=\delta^{\otimes n},\qquad Q_n=q_t^{\otimes n},
 \qquad e=(0,\ldots,0).
\]
The distinguished point has $P_n(e)=1$ and $Q_n(e)=t^n$.
The coarse-graining map
\[
 C_n(P)=\left(P(e),\sum_{\omega\ne e}P(\omega)\right)
\]
is stochastic, and sends $(P_n,Q_n)$ to $(\delta,q_{t^n})$.
Consequently
\begin{equation}\label{eq:binary:coarse}
 g(t^n)\leq D_{\C^{\Omega_n}}(P_n\relsep Q_n).
\end{equation}

For the reverse direction, define a probability distribution on the
complement of $e$ by
\[
 R_{t,n}(\omega)=
 \begin{cases}
 0,&\omega=e,\\[1mm]
 Q_n(\omega)/(1-t^n),&\omega\ne e.
 \end{cases}
\]
Its denominator is positive, and its total mass is one. The stochastic
map $B_{t,n}(a,b)=aP_n+bR_{t,n}$ is stochastic and satisfies
\[
 B_{t,n}\delta=P_n,\qquad
 B_{t,n}q_{t^n}=t^nP_n+(1-t^n)R_{t,n}=Q_n.
\]
Another application of monotonicity reverses
\eqref{eq:binary:coarse}. Combining the two directions with
\ref{ax:add}, we obtain
\begin{equation}\label{eq:binary:power}
 g(t^n)=D_{\C^{\Omega_n}}(P_n\relsep Q_n)=ng(t).
\end{equation}
The map $B_{t,n}$ depends on $t,n$, which is permitted by the axiom.
Using $g(1/2)=\log2$ in this identity gives
\begin{equation}\label{eq:binary:dyadic}
 g(2^{-m})=m\log2\quad(m\geq0),
\end{equation}
where $m=0$ uses $g(1)=0$.

\step{Step 4: determine the intermediate values.}
First $g(t)$ is finite for every $t>0$: choose $m$ with $2^{-m}\leq t$
and use monotonicity to get $g(t)\leq m\log2$.
For $0<t<1$, let $m_n=\lfloor-n\log t/\log2\rfloor$. Then
\[
 2^{-(m_n+1)}<t^n\leq2^{-m_n}.
\]
The direction of the monotonicity inequalities is therefore
\[
 g(2^{-m_n})\leq g(t^n)\leq g(2^{-(m_n+1)}).
\]
Equations~\eqref{eq:binary:power} and \eqref{eq:binary:dyadic} give
\[
 \frac{m_n\log2}{n}\leq g(t)
 \leq\frac{(m_n+1)\log2}{n}.
\]
Both bounds tend to $-\log t$. This proves the formula for $0<t<1$.
Finally $g(0)\geq g(2^{-m})=m\log2$ for every $m$, so $g(0)=+\infty$.
No continuity assumption has been used in this lemma.
\end{proof}

\subsection*{An order bound obtained from a map into \texorpdfstring{$\C^2$}{C2}}
For normal positive functionals, define
\begin{equation}\label{eq:dmax}
 \Dmax_M(\alpha\relsep\beta)
 =\inf\{\log c:c>0,\ \alpha\leq c\beta\}.
\end{equation}
The order means $\alpha(x)\leq c\beta(x)$ for every $x\geq0$.
We set $\inf\varnothing=+\infty$ and $\Dmax(0\relsep\beta)=-\infty$.
For states, evaluating an order bound at $1$ shows $c\geq1$, so
$\Dmax$ takes values in $[0,+\infty]$.

\begin{lemma}[An order bound]\label{lem:max}
For the normalized divergence, under \ref{ax:dpi} and \ref{ax:add},
\[
 D_M(\alpha\relsep\beta)\leq\Dmax_M(\alpha\relsep\beta),
 \qquad \alpha,\beta\in\Sn(M).
\]
\end{lemma}
\begin{proof}
If no finite order bound exists, the assertion has right-hand side
$+\infty$ and is automatic. Otherwise take any $c>1$ with
$\alpha\leq c\beta$. The functional
\[
 \eta=\frac{c\beta-\alpha}{c-1}
\]
is normal and positive, and $\eta(1)=(c-1)/(c-1)=1$. It is a state.
The map
\[
 \Phi:M\to\C^2,\qquad\Phi(x)=(\alpha(x),\eta(x))
\]
is normal unital completely positive: each coordinate is a normal state.
The distributions $p=(1,0)$ and $q=(c^{-1},1-c^{-1})$ satisfy
\[
 p\circ\Phi=\alpha,
 \qquad q\circ\Phi=c^{-1}\alpha+(1-c^{-1})\eta=\beta.
\]
Thus precomposition with $\Phi$ gives the desired pair exactly.
Monotonicity and Lemma~\ref{lem:binary} imply
\[
 D_M(\alpha\relsep\beta)
 \leq D_{\C^2}(p\relsep q)=-\log(c^{-1})=\log c.
\]
Taking the infimum over these $c$ proves the claim. If the infimum is
zero, one takes admissible $c>1$ decreasing to $1$; the same inequality
still gives the conclusion.
\end{proof}

\subsection*{Why lower semicontinuity becomes an asymptotic estimate}
The next lemma uses superadditivity and the fact that norm
approximation on a tensor product implies norm approximation for
the restrictions to every tensor factor. Only lower semicontinuity
on the fixed algebra $M$ is needed.

\begin{lemma}[Asymptotic lower semicontinuity]\label{lem:alsc}
Assume \ref{ax:super} and \ref{ax:lsc}. Fix $\varphi,\psi\in\Sn(M)$.
If $\omega_n\in\Sn(M^{\otimes n})$ and
$\|\omega_n-\varphi^{\otimes n}\|\to0$, then
\begin{equation}\label{eq:alsc}
 D_M(\varphi\relsep\psi)
 \leq\liminf_{n\to\infty}\frac1n
 D_{M^{\otimes n}}(\omega_n\relsep\psi^{\otimes n}).
\end{equation}
\end{lemma}
\begin{proof}
Write $\omega_{n,k}$ for the $k$th marginal and
$\varepsilon_n=\|\omega_n-\varphi^{\otimes n}\|$.
For $\|x\|\leq1$, insert $x$ into the $k$th tensor factor. The resulting
operator has norm at most one, so
\[
 |\omega_{n,k}(x)-\varphi(x)|\leq\varepsilon_n,
 \qquad
 \max_{1\leq k\leq n}\|\omega_{n,k}-\varphi\|\leq\varepsilon_n.
\]

Fix any finite real $a<D_M(\varphi\relsep\psi)$.
Lower semicontinuity in the first argument gives a radius $r>0$ such that
\[
 \|\chi-\varphi\|<r\quad\Longrightarrow\quad
 D_M(\chi\relsep\psi)>a.
\]
Indeed, if no such radius existed, states $\chi_j$ within $1/j$ of
$\varphi$ with $D_M(\chi_j\relsep\psi)\leq a$ would contradict
\ref{ax:lsc}. For large $n$, $\varepsilon_n<r$, so this estimate holds
simultaneously for every marginal $\omega_{n,k}$.
Iterating \ref{ax:super},
\[
 D_{M^{\otimes n}}(\omega_n\relsep\psi^{\otimes n})
 \geq\sum_{k=1}^nD_M(\omega_{n,k}\relsep\psi)>na.
\]
Divide by $n$, take the lower limit, and let $a$ increase to
$D_M(\varphi\relsep\psi)$. If this last value is $+\infty$, the argument
applies to every finite $a$ and proves that the lower limit is infinite.
\end{proof}

\section{The asymptotic step and its precise meaning}

This section introduces quantities defined directly from states and
their order. They do not involve the unknown function $D$.
They will bridge the order estimate in Lemma~\ref{lem:max} and the
tensor estimates in the main proof.

\subsection*{Vector representatives and approximation of functionals}
A normal positive functional $\alpha$ can be implemented by a vector
$\xi$ in a normal unital representation $\pi:M\to B(K)$, meaning
$\alpha(x)=\langle\xi,\pi(x)\xi\rangle$ for every $x\in M$.
Then $\|\xi\|^2=\alpha(1)$, so $\xi$ is a unit vector precisely when
$\alpha$ is a state.
The maximal absolute inner product of vector representatives,
usually called the root fidelity, is
\[
 F(\alpha,\beta)=\sup_{\pi,\xi,\eta}|\langle\xi,\eta\rangle|,
\]
where $\pi$ is normal and unital and $\xi,\eta$ implement the two functionals in that representation.
Both the common representation and the implementing vectors may vary.
The inner product of the two canonical natural-cone vectors need not
already be this supremum. For density matrices,
$F(\rho,\sigma)=\tr|\rho^{1/2}\sigma^{1/2}|$. 
 
If $\alpha(1),\beta(1)\leq1$, define
\begin{align*}
 \overline F(\alpha,\beta)
 &=F(\alpha,\beta)+\sqrt{(1-\alpha(1))(1-\beta(1))},\\
 P(\alpha,\beta)&=\sqrt{1-\overline F(\alpha,\beta)^2}.
\end{align*}
In all applications below, one argument is a state, so this term is zero.
The quantity $P$ is commonly called the purified distance.
For a state $\varphi$ and $0<\varepsilon<1$, put
\begin{equation}\label{eq:smooth}
 \Dmax[,\varepsilon]_M(\varphi\relsep\psi)
 =\inf\{\Dmax(\alpha\relsep\psi):
   \alpha\in M_*^+,\ \alpha(1)\leq1,\ P(\alpha,\varphi)\leq\varepsilon\}.
\end{equation}
The approximating positive functional may have value less than one at $1$.

\subsection*{Positive contractions separating states}
For $0<\varepsilon<1$, consider positive contractions $Q\in M$
on which $\varphi$ is close to one, and define
\begin{equation}\label{eq:beta}
 \beta_\varepsilon(\varphi\relsep\psi)
 =\inf\{\psi(Q):0\leq Q\leq1,\ \varphi(Q)\geq1-\varepsilon\}.
\end{equation}
Thus $\beta_\varepsilon$ is the infimum of $\psi(Q)$ subject to
$\varphi(Q)\geq1-\varepsilon$. This gives a numerical measure of
separation of the two functionals by positive contractions.

The next proposition is a known special case of the asymptotic equipartition theorem in
\cite[Theorems~3.1--3.2 and Proposition~3.6]{FGR}; the separation
estimate is due to Pautrat--Wang \cite{PW}. These results already hold
for general \vN\ algebras. Indeed, $\varphi\leq C\psi$ implies
$\Dmax(\varphi\relsep\psi)\leq\log C$, and the cited Proposition~3.6
guarantees the moment conditions in those theorems. We record the
first-order consequences in precisely the form needed here.

\begin{proposition}[Asymptotic bounds for dominated pairs]
\label{prop:aep}
Let $M$ be $\sigma$-finite and let $\varphi,\psi\in\Sn(M)$ satisfy
$\varphi\leq C\psi$ for a finite $C$. Then $s=S_M(\varphi\relsep\psi)$
is finite, and for each fixed $0<\varepsilon<1$,
\begin{align}
 \lim_{n\to\infty}\frac1n
 \Dmax[,\varepsilon]_{M^{\otimes n}}(\varphi^{\otimes n}\relsep\psi^{\otimes n})
 &=s,\label{eq:maxaep}\\
 \lim_{n\to\infty}-\frac1n
 \log\beta_\varepsilon(\varphi^{\otimes n}\relsep\psi^{\otimes n})
 &=s.\label{eq:steinaep}
\end{align}
\end{proposition}

The first limit concerns the smallest constants in inequalities
$\alpha_n\leq c_n\psi^{\otimes n}$ with $\alpha_n$ close to
$\varphi^{\otimes n}$: their logarithms, divided by $n$, tend to $s$.
The second limit describes the exponential decay of the infimum of
$\psi^{\otimes n}(Q)$ under the condition
$\varphi^{\otimes n}(Q)\geq1-\varepsilon$.
Appendix~\ref{app:aep} proves these estimates in the dominated case;
the cited papers give stronger asymptotic conclusions.


\begin{lemma}[Normalized approximants with vanishing error]
\label{lem:approximants}
Under the assumptions of Proposition~\ref{prop:aep}, there are normal
states $\omega_n$ on $M^{\otimes n}$ such that
\begin{equation}\label{eq:approximants}
 \|\omega_n-\varphi^{\otimes n}\|\to0,
 \qquad
 \limsup_{n\to\infty}\frac1n
 \Dmax(\omega_n\relsep\psi^{\otimes n})\leq s.
\end{equation}
\end{lemma}
\begin{proof}
The upper-bound construction in Appendix~\ref{app:upperorder}, using
\eqref{eq:app:hockey} and Lemma~\ref{lem:app:cutoff}, gives, for each
fixed $R>s$, normalized states $\omega_n^{(R)}$ satisfying
\[
 \|\omega_n^{(R)}-\varphi^{\otimes n}\|\to0,
 \qquad
 \limsup_n\frac1n\Dmax(\omega_n^{(R)}\relsep\psi^{\otimes n})\leq R.
\]
Use $R_j=s+1/j$. For each $j$, choose an integer $N_j$ such that,
for all $n\geq N_j$, a state $\omega_{n,j}$ from this construction obeys
\[
 \|\omega_{n,j}-\varphi^{\otimes n}\|\leq1/j,
 \qquad
 \frac1n\Dmax(\omega_{n,j}\relsep\psi^{\otimes n})\leq s+2/j.
\]
Increasing $N_j$ if needed, arrange $N_{j+1}>N_j$ and $N_j\geq j$.
For $n\geq N_1$, set
$j(n)=\max\{j:N_j\leq n\}$ and $\omega_n=\omega_{n,j(n)}$.
The set defining $j(n)$ is finite because $N_j\geq j$, and
$j(n)\to\infty$ because every $N_j$ is eventually passed.
The two estimates with $j=j(n)$ now prove
\eqref{eq:approximants}. Finitely many earlier states can be chosen
arbitrarily.
\end{proof}

\section{Detailed proof of the characterization}

Throughout the comparisons below, $D$ denotes the normalized
divergence $\widehat D$ from \eqref{eq:rescaling}. We restore the
original scale after proving equality on all normal state pairs.

\subsection{First comparison: an upper bound on dominated pairs}
For now assume that $M$ is $\sigma$-finite and $\varphi\leq C\psi$.
Let $s=S_M(\varphi\relsep\psi)<\infty$ and choose the states from
Lemma~\ref{lem:approximants}. Lemma~\ref{lem:alsc} bounds $D_M(\varphi\relsep\psi)$ by the
values of $D$ on the approximants, divided by $n$.
Lemma~\ref{lem:max} then gives the logarithmic order bound:
\begin{align*}
 D_M(\varphi\relsep\psi)
 &\leq\liminf_n\frac1n
       D_{M^{\otimes n}}(\omega_n\relsep\psi^{\otimes n})\\
 &\leq\limsup_n\frac1n
       \Dmax(\omega_n\relsep\psi^{\otimes n})
 \leq s.
\end{align*}
In particular, the unknown function is now known to be finite on
every dominated pair. No classical entropy formula for a general pair
of probability distributions has entered this comparison.

\subsection{Second comparison: positive contractions and commutative algebras}
The inequality $D_M(\varphi\relsep\psi)\geq s$ is automatic if $s=0$.
Assume $s>0$, and fix $0<R<s$. The construction of separating projections
\eqref{eq:app:directtest} in Appendix~\ref{app:aep} gives projections
$Q_n\in M^{\otimes n}$ with
\[
 h_n:=\varphi^{\otimes n}(Q_n)\to1,
 \qquad k_n:=\psi^{\otimes n}(Q_n)\leq e^{-nR}
\]
for all sufficiently large $n$. We will use only the normalized value
$D_{\C^2}(\delta\relsep u)=\log2$ to extract a lower bound from these projections.

\step{Step 1: choose a finite commutative algebra.}
Let
\[
 m_n=\left\lfloor\frac{nR}{\log2}\right\rfloor,
 \quad \Omega_n=\{0,1\}^{m_n},
 \quad a_n=\delta^{\otimes m_n},\quad u_n=u^{\otimes m_n}.
\]
Eventually $m_n\geq1$. The distribution $a_n$ is concentrated at
$e_n=0^{m_n}$, while $u_n$ assigns mass $2^{-m_n}$ to each point.
The choice of the integer ensures $k_n\leq e^{-nR}\leq2^{-m_n}$.
Therefore
\[
 r_n=\frac{u_n-k_na_n}{1-k_n}
\]
is a probability distribution: its numerator is nonnegative pointwise
and has total mass $1-k_n>0$.

\step{Step 2: define the UCP maps explicitly.}
The projection $Q_n$ defines the normal UCP map
\[
 \Theta_n:\C^2\to M^{\otimes n},\qquad
 \Theta_n(z_0,z_1)=z_0Q_n+z_1(1-Q_n).
\]
The states $a_n$ and $r_n$ define another UCP map
\[
 B_n:\C^{\Omega_n}\to\C^2,\qquad
 B_n(f)=\left(\sum_\omega a_n(\omega)f(\omega),
              \sum_\omega r_n(\omega)f(\omega)\right).
\]
Both maps are normal, unital, and completely positive. Their composition
$\Theta_n\circ B_n$ sends the second state to
\[
 k_na_n+(1-k_n)r_n=u_n,
\]
and sends the first state to
\[
 \omega_n=h_na_n+(1-h_n)r_n.
\]
Consequently
\begin{equation}\label{eq:lower:channel}
 nD_M(\varphi\relsep\psi)
 =D_{M^{\otimes n}}(\varphi^{\otimes n}\relsep\psi^{\otimes n})
 \geq D_{\C^{\Omega_n}}(\omega_n\relsep u_n).
\end{equation}
This is where tensor additivity and monotonicity are used.

\step{Step 3: compare all binary marginals with the normalized pair.}
Since $\omega_n-a_n=(1-h_n)(r_n-a_n)$ and the distance between two
probability distributions is at most $2$,
\[
 \|\omega_n-a_n\|_1\leq2(1-h_n)\to0.
\]
Every marginal $\omega_{n,k}$ is therefore within $2(1-h_n)$ of
$\delta$. Fix $\eta>0$. Lower semicontinuity in the first argument at
the fixed pair
$(\delta,u)$ gives a radius $r>0$ such that
\[
 \|\chi-\delta\|_1<r
 \quad\Longrightarrow\quad
 D_{\C^2}(\chi\relsep u)>\log2-\eta.
\]
For large $n$, this applies simultaneously to all $m_n$ marginals.
Iterating superadditivity in \eqref{eq:lower:channel},
\[
 nD_M(\varphi\relsep\psi)
 \geq\sum_{k=1}^{m_n}D_{\C^2}(\omega_{n,k}\relsep u)
 \geq m_n(\log2-\eta).
\]
Divide by $n$ and use $m_n/n\to R/\log2$. Then let $\eta\downarrow0$
and $R\uparrow s$. We obtain $D_M(\varphi\relsep\psi)\geq s$.
Together with the upper comparison, this proves equality on dominated
pairs in $\sigma$-finite algebras.

\subsection{Remove domination by convex combinations}
Keep $M$ $\sigma$-finite, but let $\varphi,\psi$ now be arbitrary.
For $0<t<1$, put
\[
 \psi_t=(1-t)\psi+t\varphi.
\]
The order inequality $\varphi\leq t^{-1}\psi_t$ shows that the result
already proved applies to $(\varphi,\psi_t)$:
\begin{equation}\label{eq:regularized}
 D_M(\varphi\relsep\psi_t)=S_M(\varphi\relsep\psi_t).
\end{equation}
There are two inequalities to pass to the limit, and they come from
different axioms.

The map
\[
 T_t:M\to M,\qquad T_t(x)=(1-t)x+t\varphi(x)1
\]
is normal, unital, and completely positive. Indeed, it is a convex combination
of the identity and state evaluation followed by scalar inclusion.
It satisfies
\[
 \varphi\circ T_t=\varphi,
 \qquad \psi\circ T_t=(1-t)\psi+t\varphi=\psi_t.
\]
Monotonicity gives the lower estimate
\begin{equation}\label{eq:regularizedlower}
 D_M(\varphi\relsep\psi)\geq S_M(\varphi\relsep\psi_t).
\end{equation}
For the upper estimate, note that
$\|\psi_t-\psi\|=t\|\varphi-\psi\|\leq2t$.
Lower semicontinuity in the second argument, with $\varphi$ fixed,
along any sequence $t\downarrow0$ gives
\begin{equation}\label{eq:regularizedupper}
 D_M(\varphi\relsep\psi)
 \leq\liminf_{t\downarrow0}S_M(\varphi\relsep\psi_t).
\end{equation}

We now use the basic convexity and lower semicontinuity properties of
Araki relative entropy to evaluate this limit. Lower semicontinuity
always gives
$S_M(\varphi\relsep\psi)\leq\liminf_{t\downarrow0}S_M(\varphi\relsep\psi_t)$.
If $S_M(\varphi\relsep\psi)<\infty$, convexity in the second argument gives
\[
 S_M(\varphi\relsep\psi_t)
 \leq(1-t)S_M(\varphi\relsep\psi)+tS_M(\varphi\relsep\varphi)
 =(1-t)S_M(\varphi\relsep\psi).
\]
If $S_M(\varphi\relsep\psi)=+\infty$, the lower-semicontinuity inequality
already forces the lower limit to be $+\infty$.
Thus in both cases
\begin{equation}\label{eq:entropyreg}
 \lim_{t\downarrow0}S_M(\varphi\relsep\psi_t)=S_M(\varphi\relsep\psi).
\end{equation}
Combining \eqref{eq:regularizedlower}, \eqref{eq:regularizedupper}, and
\eqref{eq:entropyreg} proves equality for every normal pair on $M$.
All comparisons use nonnegative extended quantities; no subtraction
of infinities is involved.

\subsection{Remove \texorpdfstring{$\sigma$}{sigma}-finiteness by a support corner}
Let $M$ be arbitrary, and put
\[
 p=\supp(\varphi)\vee\supp(\psi),\qquad M_0=pMp.
\]
The restrictions $\varphi_0,\psi_0$ are states on $M_0$, whose unit is
$p$. The state $(\varphi_0+\psi_0)/2$ has support $p$ and is faithful
on this corner. Hence $M_0$ is $\sigma$-finite.

Because both original states are supported on $p$,
$\varphi(x)=\varphi(pxp)$ and $\psi(x)=\psi(pxp)$.
For example, the mixed terms involving $1-p$ vanish by
Cauchy--Schwarz for a positive functional. Consider compression
\[
 C:M\to M_0,\qquad C(x)=pxp.
\]
It is normal unital completely positive, with $C(1)=p$, and
$\varphi=\varphi_0\circ C$, $\psi=\psi_0\circ C$.
Choose a normal state $\theta$ on $M_0$, for instance $\varphi_0$, and set
\[
 J:M_0\to M,\qquad J(y)=y+\theta(y)(1-p).
\]
This map is normal and completely positive, as the block diagonal sum
of inclusion on $p$ and the scalar map determined by $\theta$ on $1-p$.
It is unital because $J(p)=p+(1-p)=1$, and $C\circ J=\id_{M_0}$.
Consequently $\varphi\circ J=\varphi_0$ and $\psi\circ J=\psi_0$.
The two instances of monotonicity now give
\[
 D_M(\varphi\relsep\psi)
 \leq D_{M_0}(\varphi_0\relsep\psi_0)
 \leq D_M(\varphi\relsep\psi).
\]
Araki relative entropy obeys the same two inequalities by its own
monotonicity property. Equality on $M_0$ therefore proves equality on
$M$. We have proved $\widehat D=S$ on every normal state pair.
Restoring the original divergence in \eqref{eq:rescaling} gives
\[
 D_M(\varphi\relsep\psi)=\frac{b}{\log2}S_M(\varphi\relsep\psi),
\]
where $0<b<\infty$ was established in Lemma~\ref{lem:proper}.

\subsection{Why Araki relative entropy satisfies the axioms}
Direct evaluation gives $S_{\C^2}(\delta\relsep u)=\log2$:
the only occupied coordinate contributes $\log2$. Thus $S$ is non-degenerate.
Monotonicity and separate norm lower semicontinuity are standard
theorems for Araki relative entropy \cite{Araki76,Araki77,OP}.
Tensor additivity follows from the tensor-product identity for relative
modular operators: the logarithm of the tensor product is the sum of
the two logarithms on their supports.

For completeness, the specific superadditivity needed here follows
from the conditional-expectation chain rule for a product second state:
\begin{align*}
 S_{M\ot N}(\omega\relsep\psi_M\otimes\psi_N)
 &=S_M(\omega_M\relsep\psi_M)
   +S_{M\ot N}(\omega\relsep\omega_M\otimes\psi_N)\\
 &\geq S_M(\omega_M\relsep\psi_M)+S_N(\omega_N\relsep\psi_N).
\end{align*}
The relevant expectation is the slice map $\id_M\otimes\psi_N$ onto
$M\otimes1$. Its extension of $\omega_M$ is
$\omega_M\otimes\psi_N$, which explains the middle second state.
The last inequality is monotonicity under restriction to $N$.
The chain rule is an equality of extended nonnegative quantities, so
this verification also covers infinite values. Multiplication by a
positive finite constant preserves non-degeneracy and all four axioms.
This completes the proof of Theorem~\ref{thm:main}.
\qed

\section{Alternative axioms: conditional expectations}
\label{sec:alternative}

\subsection{The conditional-expectation chain rule}
Let $P\subseteq R$ be a unital inclusion of \vN\ algebras. A normal
conditional expectation $E:R\to P$ is a normal unital completely
positive projection onto $P$ satisfying
\[
 E(a x b)=aE(x)b,\qquad a,b\in P,\ x\in R.
\]
For a state $\omega$ on $R$ write $\omega_P=\omega|_P$.
Its $E$-invariant extension back to $R$ using $E$ is $\omega_P\circ E$.

Consider the following axiom, imposed whenever such a normal expectation
exists:
\begin{quote}
\textbf{(CE) Conditional-expectation chain rule.}
If $\psi=\psi_P\circ E$, then, for every normal state $\omega$,
\begin{equation}\label{eq:ce:general}
 D_R(\omega\relsep\psi)
 =D_P(\omega_P\relsep\psi_P)
  +D_R(\omega\relsep\omega_P\circ E).
\end{equation}
\end{quote}
This is an equality of nonnegative extended quantities. Neither the
states nor the expectation are required to be faithful. This full
domain is part of the formulation, including when a term is infinite.

The conditional-expectation property is a classical ingredient of
Petz's axiomatic approach; see \cite{OP,Petz92,Petz2006}.
Here it gives an alternative to the two tensor axioms.

\begin{theorem}[Alternative characterization]\label{thm:alternative}
Let $D$ be a non-degenerate divergence on all normal state pairs of all
\vN\ algebras. The following three assumptions imply
$D=\lambda S$ for a universal constant $0<\lambda<\infty$:
\begin{enumerate}[label=\textup{(C\arabic*)},ref=C\arabic*]
 \item monotonicity \ref{ax:dpi};
 \item the chain rule \textup{(CE)};
 \item separate norm lower semicontinuity \ref{ax:lsc}.
\end{enumerate}
Conversely, every positive finite multiple of Araki relative entropy
is non-degenerate and satisfies these three assumptions.
For uniqueness, \textup{(CE)} can even be restricted to the tensor-factor
expectations used below.
\end{theorem}
\begin{proof}
We derive \ref{ax:add} and \ref{ax:super}, then apply
Theorem~\ref{thm:main}. The argument uses no additivity assumption
before it has been proved.

\step{Step 1: adjoining a common state preserves the value of $D$.}
Let $\chi\in\Sn(M)$ and $\alpha,\beta\in\Sn(N)$. The maps
\[
 \iota:N\to M\ot N,\quad\iota(y)=1\otimes y,
 \qquad
 P_\chi:M\ot N\to N,\quad P_\chi=\chi\otimes\id_N
\]
are normal unital completely positive, and $P_\chi\circ\iota=\id_N$.
Restriction along $\iota$ sends
$(\chi\otimes\alpha,\chi\otimes\beta)$ to $(\alpha,\beta)$;
composition with $P_\chi$ sends the latter pair back. Thus monotonicity in both directions gives
\begin{equation}\label{eq:ce:commonfactor}
 D_{M\ot N}(\chi\otimes\alpha\relsep\chi\otimes\beta)
 =D_N(\alpha\relsep\beta).
\end{equation}
This conclusion holds even if the common value is infinite.

\step{Step 2: apply the chain rule to a tensor-factor expectation.}
For fixed $\psi_N\in\Sn(N)$ let
\[
 E_{\psi_N}:M\ot N\to M\otimes1,\qquad
 E_{\psi_N}(x\otimes y)=\psi_N(y)\,x\otimes1.
\]
This is a normal conditional expectation. Its range is identified with
$M$ by $x\mapsto x\otimes1$. This identification preserves $D$ by
monotonicity for an isomorphism and its inverse.
The second state $\psi_M\otimes\psi_N$ is invariant under $E_{\psi_N}$,
whereas $\omega_M\circ E_{\psi_N}=\omega_M\otimes\psi_N$.
Consequently \eqref{eq:ce:general} reads
\begin{equation}\label{eq:ce:product}
 \begin{split}
 D_{M\ot N}(\omega\relsep\psi_M\otimes\psi_N)
 ={}&D_M(\omega_M\relsep\psi_M)\\
 &+D_{M\ot N}(\omega\relsep\omega_M\otimes\psi_N).
 \end{split}
\end{equation}

\step{Step 3: obtain tensor additivity.}
Take $\omega=\varphi_M\otimes\varphi_N$ in
\eqref{eq:ce:product}. Equation~\eqref{eq:ce:commonfactor} evaluates
the last term as $D_N(\varphi_N\relsep\psi_N)$, so
\[
 D_{M\ot N}(\varphi_M\otimes\varphi_N\relsep\psi_M\otimes\psi_N)
 =D_M(\varphi_M\relsep\psi_M)+D_N(\varphi_N\relsep\psi_N).
\]
This is \ref{ax:add}.

\step{Step 4: obtain superadditivity.}
For arbitrary $\omega$, restrict the last term of
\eqref{eq:ce:product} to $1\otimes N$. Monotonicity gives
\[
 D_{M\ot N}(\omega\relsep\omega_M\otimes\psi_N)
 \geq D_N(\omega_N\relsep\psi_N).
\]
Together with \eqref{eq:ce:product}, this is \ref{ax:super}.
The original non-degenerate divergence now satisfies all four axioms of
Theorem~\ref{thm:main}, which gives $D=\lambda S$, with
$\lambda=b/\log2$. Conversely, Araki entropy has the
conditional-expectation decomposition \eqref{eq:ce:general}
\cite{Araki77,OP,Petz2006}, is non-degenerate, and has the other two
properties. These assertions persist under multiplication by any
positive finite constant.
\end{proof}

The proof only uses \eqref{eq:ce:product}, for all algebras and all
normal state choices. Thus that restricted identity can replace the
full axiom \textup{(CE)} in Theorem~\ref{thm:alternative}. This keeps
the alternative assumption close to the product structure of the second state
actually needed in the proof.

\section{States of C*-algebras}\label{sec:cstar}

\subsection{The entropy and the C*-algebra axioms}
Let $A$ be a unital C*-algebra and write $\mathcal S(A)$ for its
state space. Every state $\varphi$ has a unique normal extension
$\widetilde\varphi$ to the enveloping \vN\ algebra $A^{**}$. The
correspondence $A^*=(A^{**})_*$ is isometric, so in particular
\[
 \|\widetilde\varphi-\widetilde\psi\|=\|\varphi-\psi\|.
\]
Define the standard relative entropy by
\begin{equation}\label{eq:cstar:entropy}
 S_A(\varphi\relsep\psi)
 :=S_{A^{**}}(\widetilde\varphi\relsep\widetilde\psi).
\end{equation}
Equivalently, if both states have normal extensions to
$\pi(A)''$ in a common nondegenerate representation, their entropy
can be computed there. This compatibility with representations, and
the agreement with Uhlmann's construction, are classical; see
\cite[Section~3, especially Lemma~3.1]{HOT83}.

Throughout this section $A\otimes_{\min}B$ denotes the minimal
C*-tensor product. All its states and all unital completely
positive maps between unital C*-algebras are allowed. There
is no normality qualification in this formulation.
Non-degeneracy has the same meaning as before: there are a unital
C*-algebra $A_0$ and states $\varphi_0,\psi_0$ with
$0<D_{A_0}(\varphi_0\relsep\psi_0)<\infty$.

\begin{theorem}[C*-algebra characterization]\label{thm:cstar}
Let $D=(D_A)_A$ be a family
\[
 D_A:\mathcal S(A)\times\mathcal S(A)\longrightarrow[0,+\infty]
\]
indexed by all unital C*-algebras. Assume that $D$ is non-degenerate and
satisfies the following four conditions.
\begin{enumerate}[label=\textup{(B\arabic*)},ref=B\arabic*]
\item\label{ax:cstar:dpi}\textbf{Monotonicity.}
For every UCP map $\Phi:B\to A$ and states $\varphi,\psi$ of $A$,
\[
 D_B(\varphi\circ\Phi\relsep\psi\circ\Phi)
 \leq D_A(\varphi\relsep\psi).
\]
\item\label{ax:cstar:add}\textbf{Tensor additivity.}
For states $\varphi_i,\psi_i$ of $A_i$,
\[
 D_{A_1\otimes_{\min}A_2}
 (\varphi_1\otimes\varphi_2\relsep\psi_1\otimes\psi_2)
 =D_{A_1}(\varphi_1\relsep\psi_1)+D_{A_2}(\varphi_2\relsep\psi_2).
\]
\item\label{ax:cstar:super}\textbf{Superadditivity for product second states.}
For every state $\omega$ of $A\otimes_{\min}B$,
\[
 D_{A\otimes_{\min}B}(\omega\relsep\psi_A\otimes\psi_B)
 \geq D_A(\omega_A\relsep\psi_A)+D_B(\omega_B\relsep\psi_B).
\]
Here the marginals are given by restriction to $A\otimes1$ and
$1\otimes B$.
\item\label{ax:cstar:lsc}\textbf{Separate norm lower semicontinuity.}
With either state fixed, $D_A$ is lower semicontinuous in the other
state for the norm of $A^*$.
\end{enumerate}
Then there is a universal constant $0<\lambda<\infty$ such that
\[
 D_A(\varphi\relsep\psi)=\lambda S_A(\varphi\relsep\psi)
\]
for every state pair on every unital C*-algebra. It is given by
$\lambda=b/\log2$, where $0<b:=D_{\C^2}(\delta\relsep u)<\infty$.
Conversely, every positive finite multiple of $S$ is non-degenerate and
satisfies these conditions. No nuclearity or faithfulness assumption
is required.
\end{theorem}

As in Remark~\ref{rem:scale}, prescribing the optional binary value
$\log2$ selects $D=S$. Weak additivity on identical
tensor powers is enough for the uniqueness proof. We give the
transfer argument explicitly, because the asymptotic objects first
arise in a \vN\ algebra whereas the unknown function is evaluated
on the original C*-algebras.

\subsection{Restriction of states and approximation of positive contractions}
The following lemma is elementary.  

\begin{lemma}[Density and restriction]\label{lem:cstar:density}
Put $M=A^{**}$, $A_n=A^{\otimes_{\min}n}$, and
$N_n=M^{\otimes n}$, with the latter tensor product spatial.
There is a canonical injective unital *-homomorphism
\[
 j_n:A_n\longrightarrow N_n
\]
whose range is strongly dense in a faithful normal representation
of $N_n$. The normal product states on $N_n$ restrict to the
corresponding product states on $A_n$.

Restriction of normal positive functionals along $j_n$ preserves
order inequalities and does not increase norm distances. Moreover,
given a positive contraction $0\leq Q\leq1$ in $N_n$, normal states $\alpha,\beta$
of $N_n$, and $\varepsilon>0$, there is a positive contraction $a\in A_n$ such that
\begin{equation}\label{eq:cstar:effectapprox}
 |\alpha(j_n(a))-\alpha(Q)|<\varepsilon,
 \qquad
 |\beta(j_n(a))-\beta(Q)|<\varepsilon.
\end{equation}
\end{lemma}
\begin{proof}
Represent $M$ faithfully and normally on a Hilbert space $H$.
Its restriction to $A$ is faithful. The representation of the
algebraic tensor product $A^{\odot n}$ on $H^{\otimes n}$ therefore
defines its minimal C*-norm and extends to the injective map $j_n$.
Since $A$ is strongly dense in $M$, the range of $j_n$ generates
$M^{\otimes n}$ as a \vN\ algebra. Evaluation on elementary tensors
identifies the restricted product states.

If $\varphi\leq c\psi$ on $N_n$, evaluation at $j_n(x)$ for $x\geq0$
gives the same order bound after restriction. The norm estimate
follows because $j_n$ maps the unit ball into the unit ball.
Finally Kaplansky density makes the positive unit ball of $j_n(A_n)$
strongly dense in the positive unit ball of $N_n$. A bounded net
converging strongly also converges ultraweakly, so its values under
the two normal states converge. Choose one member of that net to
obtain \eqref{eq:cstar:effectapprox}. No metrizability or separability
assumption is needed for this finite collection of evaluations.
\end{proof}

For states on a C*-algebra we use the order definition
\[
 \Dmax_A(\varphi\relsep\psi)
 =\inf\{\log c:c>0,\ \varphi\leq c\psi\}.
\]
The next lemma records exactly the two consequences of the known
\vN\ asymptotic results that enter the proof.

\begin{lemma}[Dominated C*-algebra pairs]\label{lem:cstar:asymptotic}
Let $\varphi,\psi\in\mathcal S(A)$ satisfy $\varphi\leq C\psi$ for
a finite $C$, and put $s=S_A(\varphi\relsep\psi)$. Then
$0\leq s\leq\log C$, and there are states $\omega_n$ on $A_n$ with
\begin{equation}\label{eq:cstar:approximants}
 \|\omega_n-\varphi^{\otimes n}\|\to0,
 \qquad
 \limsup_n\frac1n
 \Dmax_{A_n}(\omega_n\relsep\psi^{\otimes n})\leq s.
\end{equation}
For every $0<R<s$ there are positive contractions $a_n\in A_n$ such that
\begin{equation}\label{eq:cstar:tests}
 \varphi^{\otimes n}(a_n)\to1,
 \qquad \psi^{\otimes n}(a_n)\leq2e^{-nR}
\end{equation}
for all sufficiently large $n$.
\end{lemma}
\begin{proof}
The order bound extends to $M=A^{**}$. Let
$p=\supp(\widetilde\psi)$. Domination implies
$\supp(\widetilde\varphi)\leq p$, so both normal states are supported
on the corner $M_0=pMp$, where $\widetilde\psi$ is faithful and
$M_0$ is $\sigma$-finite. Entropy on this corner equals
$S_A(\varphi\relsep\psi)$ by the support-corner argument of Section~4.
Proposition~\ref{prop:aep} and Lemma~\ref{lem:approximants} apply on
$M_0$. They give the finite bound on $s$ and normal approximating
states on $M_0^{\otimes n}$.

Write $p_n=p^{\otimes n}$. Extend these approximating
states normally to $N_n$ by compression
$x\mapsto p_n x p_n$.
Compression preserves their state normalization, their order bounds
against $\widetilde\psi^{\otimes n}$, and their norm approximation to
$\widetilde\varphi^{\otimes n}$. Restrict next to $A_n$ along $j_n$.
Lemma~\ref{lem:cstar:density} gives both assertions of
\eqref{eq:cstar:approximants}.

For the separating projections, the construction
\eqref{eq:app:directtest} on $M_0$ gives projections $Q_n$ with
\[
 \widetilde\varphi^{\otimes n}(Q_n)\to1,
 \qquad \widetilde\psi^{\otimes n}(Q_n)\leq e^{-nR}
\]
eventually. Regard $Q_n$ as positive contractions of $N_n$ supported on
$p_n$. Apply Lemma~\ref{lem:cstar:density} separately
for each $n$, with accuracy $e^{-2nR}$ in the two state evaluations.
The resulting positive contractions in $A_n$ satisfy
$\varphi^{\otimes n}(a_n)\to1$ and
$\psi^{\otimes n}(a_n)$ is at most
$e^{-nR}+e^{-2nR}\leq2e^{-nR}$, proving
\eqref{eq:cstar:tests}. The exponentially small choice of accuracy
preserves the exponential estimate.
\end{proof}

This lemma transfers established asymptotic estimates to the
C*-algebra setting. The input is the dominated case of
\cite{FGR,PW}, reproved in Appendix~\ref{app:aep}. The tensor products
are handled through the dense inclusion $j_n$. 

\subsection{Proof of the C*-algebra characterization}
\begin{proof}[Proof of Theorem~\ref{thm:cstar}]
\step{Step 1: the elementary consequences of the axioms.}
The proof of Lemma~\ref{lem:proper} uses only states, UCP maps,
tensor powers, and separate norm lower semicontinuity. All these
arguments are valid for C*-algebras, so non-degeneracy and
\ref{ax:cstar:dpi}--\ref{ax:cstar:lsc} imply $0<b<\infty$.
For Steps~1--4 of this proof replace $D$ by $(\log2/b)D$ and retain
the notation $D$ for this normalized divergence.
The proofs of Lemmas~\ref{lem:binary} and \ref{lem:max} also use only
states, UCP maps, and tensor powers. Thus their conclusions hold
under \ref{ax:cstar:dpi} and \ref{ax:cstar:add} after this normalization.
In particular,
\[
 D_A(\chi\relsep\chi)=0,\qquad
 D_{\C^2}(\delta\relsep(t,1-t))=-\log t,
 \qquad D_A(\rho\relsep\psi)\leq\Dmax_A(\rho\relsep\psi).
\]
The norm estimate for marginals in Lemma~\ref{lem:alsc} also uses only
contractivity under restriction. Thus
\ref{ax:cstar:super} and \ref{ax:cstar:lsc} imply its C*-algebra version:
if $\omega_n\in\mathcal S(A_n)$ and
$\|\omega_n-\varphi^{\otimes n}\|\to0$, then
\begin{equation}\label{eq:cstar:alsc}
 D_A(\varphi\relsep\psi)
 \leq\liminf_n\frac1n
 D_{A_n}(\omega_n\relsep\psi^{\otimes n}).
\end{equation}

\step{Step 2: the upper comparison for dominated pairs.}
Assume $\varphi\leq C\psi$ and use the states from
Lemma~\ref{lem:cstar:asymptotic}. Combining
\eqref{eq:cstar:alsc} with the order bound gives
\[
 D_A(\varphi\relsep\psi)
 \leq\liminf_n\frac1nD_{A_n}(\omega_n\relsep\psi^{\otimes n})
 \leq\limsup_n\frac1n\Dmax_{A_n}(\omega_n\relsep\psi^{\otimes n})
 \leq s.
\]
Thus the unknown function is finite on these pairs.

\step{Step 3: the lower comparison for dominated pairs.}
If $s=0$, nonnegativity gives the lower comparison. Otherwise fix
$0<R<s$ and choose the positive contractions in \eqref{eq:cstar:tests}.
Write $h_n=\varphi^{\otimes n}(a_n)$ and
$k_n=\psi^{\otimes n}(a_n)$, and set
\[
 m_n=\max\left\{0,
       \left\lfloor\frac{nR-\log2}{\log2}\right\rfloor\right\}.
\]
For all large $n$, $m_n\geq1$, $m_n/n\to R/\log2$, and
$k_n\leq2e^{-nR}\leq2^{-m_n}$. On $\{0,1\}^{m_n}$ put
$d_n=\delta^{\otimes m_n}$, $u_n=u^{\otimes m_n}$, and
\[
 r_n=\frac{u_n-k_nd_n}{1-k_n},
 \qquad \nu_n=h_nd_n+(1-h_n)r_n.
\]
The first expression is a probability distribution because
$k_n\leq2^{-m_n}$. The UCP maps determined by $a_n$ and by the two states $d_n,r_n$
carry the pair
$(\varphi^{\otimes n},\psi^{\otimes n})$ to $(\nu_n,u_n)$.
The UCP map is the composition of
\[
 (z_0,z_1)\longmapsto z_0a_n+z_1(1-a_n)
\]
with the same map between finite commutative algebras as in Section~4. Both maps
are UCP, even though $a_n$ need not be a projection. Monotonicity
and tensor additivity give
\[
 nD_A(\varphi\relsep\psi)
 \geq D_{\C^{2^{m_n}}}(\nu_n\relsep u_n).
\]
Since $\|\nu_n-d_n\|_1\leq2(1-h_n)\to0$, each binary marginal of
$\nu_n$ tends uniformly to $\delta$. Lower semicontinuity at
$(\delta,u)$ and iterated superadditivity imply, for every $\eta>0$
and all large $n$,
\[
 nD_A(\varphi\relsep\psi)\geq m_n(\log2-\eta).
\]
Divide by $n$, let $\eta\downarrow0$, and then let $R\uparrow s$.
This proves $D_A(\varphi\relsep\psi)\geq s$, hence equality for all
dominated pairs.

\step{Step 4: remove domination.}
For arbitrary states set $\psi_t=(1-t)\psi+t\varphi$, $0<t<1$.
Since $\varphi\leq t^{-1}\psi_t$, the preceding steps give
$D_A(\varphi\relsep\psi_t)=S_A(\varphi\relsep\psi_t)$.
The UCP map $T_t:A\to A$, $T_t(x)=(1-t)x+t\varphi(x)1$, fixes
$\varphi$ and sends $\psi$ to $\psi_t$ under composition. Therefore
\[
 D_A(\varphi\relsep\psi)\geq S_A(\varphi\relsep\psi_t).
\]
Because $\|\psi_t-\psi\|\leq2t$, separate norm lower
semicontinuity gives the reverse limiting estimate
\[
 D_A(\varphi\relsep\psi)
 \leq\liminf_{t\downarrow0}S_A(\varphi\relsep\psi_t).
\]
The normal extensions satisfy
$\widetilde\psi_t=(1-t)\widetilde\psi+t\widetilde\varphi$.
The convexity and norm lower semicontinuity of Araki entropy on
$A^{**}$ give, exactly as in \eqref{eq:entropyreg},
\[
 \lim_{t\downarrow0}S_A(\varphi\relsep\psi_t)=S_A(\varphi\relsep\psi).
\]
Both inequalities therefore give equality for the normalized
divergence on every state pair, including those with infinite entropy.
Undoing the normalization gives $D_A=(b/\log2)S_A$ for the original
family. Only first-argument norm lower
semicontinuity and second-argument lower semicontinuity along these
mixtures were used.

\step{Step 5: verify the converse.}
For a UCP map $\Phi:B\to A$, the bidual map
$\Phi^{**}:B^{**}\to A^{**}$ is normal UCP, and
$\widetilde{\varphi\circ\Phi}=\widetilde\varphi\circ\Phi^{**}$.
Normal monotonicity thus gives \ref{ax:cstar:dpi} for $S$.
The scalar computation $S_{\C^2}(\delta\relsep u)=\log2$ proves
non-degeneracy, and the isometry
$A^*=(A^{**})_*$ transfers norm lower semicontinuity, giving
\ref{ax:cstar:lsc}.

For the tensor properties, it is convenient first to verify the
C*-algebra chain rule. If $B\subseteq A$ is a unital C*-subalgebra
and $E:A\to B$ is a conditional expectation, identify $B^{**}$ with
its canonical image in $A^{**}$. Then $E^{**}$ is a normal conditional
expectation onto $B^{**}$. If $\psi=\psi|_B\circ E$, the normal
extensions have the same invariance. The normal chain rule of
Section~\ref{sec:alternative} consequently gives
\begin{equation}\label{eq:cstar:ce}
 S_A(\omega\relsep\psi)
 =S_B(\omega|_B\relsep\psi|_B)
  +S_A\bigl(\omega\relsep(\omega|_B)\circ E\bigr).
\end{equation}
This familiar C*-algebra conditional-expectation property is also
listed in \cite{Petz1999,Petz2006}. It is an equality of extended
nonnegative quantities, with no faithfulness restriction.
For $A\otimes_{\min}B$ take the slice expectation
\[
 E_{\psi_B}(a\otimes b)=\psi_B(b)a\otimes1.
\]
The common-factor identity follows from adjoining a fixed state and
restricting back, using monotonicity in both directions. Hence
\eqref{eq:cstar:ce} gives tensor additivity on product first states,
and gives superadditivity on arbitrary first states by restriction
of its last term to $B$. This is the argument of
Theorem~\ref{thm:alternative}, now with minimal C*-tensor products.
It proves \ref{ax:cstar:add} and \ref{ax:cstar:super}. Every positive
finite multiple preserves these four properties and non-degeneracy,
completing the converse.
\end{proof}

\Needspace{18\baselineskip}
\subsection{The conditional-expectation formulation}
\begin{theorem}[Alternative C*-algebra characterization]
\label{thm:cstar:alternative}
For a non-degenerate divergence $D$ on all state pairs of all unital
C*-algebras, the following three conditions imply $D=\lambda S$
for one universal constant $0<\lambda<\infty$:
\begin{enumerate}
\item monotonicity \ref{ax:cstar:dpi};
\item the chain rule
\[
 D_A(\omega\relsep\psi)
 =D_B(\omega|_B\relsep\psi|_B)
  +D_A\bigl(\omega\relsep(\omega|_B)\circ E\bigr)
\]
for every conditional expectation $E:A\to B$ and every
$\psi=\psi|_B\circ E$;
\item separate norm lower semicontinuity \ref{ax:cstar:lsc}.
\end{enumerate}
Conversely, every positive finite multiple of $S$ is non-degenerate and
satisfies all three. For uniqueness the chain rule is needed only for
the tensor-factor slice expectations. The constant is again
$\lambda=b/\log2$.
\end{theorem}
\begin{proof}
The proof of Theorem~\ref{thm:alternative} is valid for C*-algebras:
the factor inclusions and the slice maps are UCP for the minimal
tensor product. Applying the chain rule to
$\id_A\otimes\psi_B$ gives tensor additivity and superadditivity.
Theorem~\ref{thm:cstar} then gives uniqueness. The converse follows
from \eqref{eq:cstar:ce} and the other properties just proved.
\end{proof}

\begin{remark}[Non-unital C*-algebras]\label{rem:cstar:nonunital}
The characterization in Theorem~\ref{thm:cstar} also holds for arbitrary,
not necessarily unital, C*-algebras. States are positive functionals
of norm one. The monotonicity axiom is reformulated for completely
positive maps $\Phi:B\to A$ whose pullbacks preserve states, that is,
$\rho\circ\Phi\in\mathcal S(B)$ for every $\rho\in\mathcal S(A)$:
\[
 D_B(\varphi\circ\Phi\relsep\psi\circ\Phi)
 \leq D_A(\varphi\relsep\psi).
\]
The minimal tensor product is retained, with marginals defined using
the canonical extensions of states to the multiplier algebras.
Together with non-degeneracy and the other three axioms, this again
characterizes $D_A=\lambda S_A$ for a universal constant
$0<\lambda<\infty$, where $S_A$ is defined by
\eqref{eq:cstar:entropy}. We omit the proof of this extension.
\end{remark}

\Needspace{8\baselineskip}
\section*{Acknowledgements}
During the preparation of this note, the author used ChatGPT (OpenAI)
for mathematical discussion, assistance with developing and drafting
proofs, identifying references, and preparing the \LaTeX{} manuscript.
The author takes full responsibility for the mathematical content,
references, and final text.

\appendix
\numberwithin{equation}{section}
\section{A detailed proof of Proposition~\ref{prop:aep}}\label{app:aep}

For the convenience of the reader, 
this appendix gives a direct, self-contained proof of the known special
case stated in Proposition~\ref{prop:aep}. Its purpose is to derive
precisely the analytic inputs needed in Section~4 from standard modular
and $L^p$ facts. The proposition already follows from
\cite[Theorems~3.1--3.2 and Proposition~3.6]{FGR}, with the separation
estimates due to \cite{PW}. We make no priority claim for the proof
techniques assembled here.

The route is tailored to domination and first-order limits: a Chernoff
bound gives separating projections, and a bound on a positive part gives
a nearby normalized state with an order bound. It uses no second-order
Gaussian expansion. The positive-power proof of the Chernoff bound is
Ozawa's argument, printed in \cite[Proposition~1.1]{JOPS}; the general
\vN\ algebra inequality and related estimates for positive contractions are
developed in \cite{JOPS,Ogata2010}.

Throughout this appendix $M$ is $\sigma$-finite, $\varphi,\psi$ are
normal states, and
\[
 \varphi\leq C\psi,\qquad 1\leq C<\infty.
\]
We write $\varphi_n=\varphi^{\otimes n}$ and
$\psi_n=\psi^{\otimes n}$. The order bound is stronger than mere support
inclusion. Its role is to control both tails of the relative modular
logarithm, so that a first-order expansion of its exponential moment is
available.

\subsection{The standard-form and \texorpdfstring{$L^p$}{Lp} ingredients}
We recall the background facts used below. In Haagerup's construction,
a normal positive functional $\chi$ corresponds to a positive density
$h_\chi\in L^1(M)$. The correspondence preserves order, addition, and
the Jordan decomposition. Its canonical trace pairing satisfies
\[
 \htr(h_\chi)=\chi(1),\qquad
 \chi(x)=\htr(h_\chi x),\qquad
 \xi_\chi=h_\chi^{1/2}\in L^2(M).
\]
The space $L^2(M)$ realizes the standard form of $M$.
The symbol $\htr$ is the canonical pairing on these $L^p$ spaces; it
does not assert that $M$ itself has a trace. We use H\"older's
inequality, the cyclicity of this pairing when the products are in
$L^1$, and positivity of the pairing of positive dual elements. See
Haagerup \cite{HaagerupLp} and Hiai \cite[Appendix~A]{Hiai2018}.

As is well known, for $0<t<1$, positive powers preserve order. One way to see the
operator-theoretic reason is the integral representation
\[
 A^t=\frac{\sin\pi t}{\pi}\int_0^\infty
 \lambda^{t-1}A(A+\lambda)^{-1}\,d\lambda.
\]
The resolvent reverses order, so $A(A+\lambda)^{-1}
=1-\lambda(A+\lambda)^{-1}$ preserves it. The formula and its order
consequence extend to the positive affiliated operators used here by
the usual quadratic-form interpretation.

The link between these densities and relative modular operators is
\begin{equation}\label{eq:app:power}
 \htr(h_\chi^{1-t}h_\eta^t)
 =\langle\xi_\chi,\Delta_{\eta,\chi}^{\,t}\xi_\chi\rangle,
 \qquad 0<t<1.
\end{equation}
In the matrix case this is just the identity obtained from
$\Delta_{\eta,\chi}=L_{h_\eta}R_{h_\chi}^{-1}$ on the appropriate
support. These facts are the standard-form and $L^p$ background; the
estimates and limiting arguments needed for the proposition are proved
below.

\subsection{An inequality for positive contractions and a vector estimate}
\begin{lemma}[Chernoff inequality]\label{lem:app:chernoff}
For normal positive functionals $\chi,\eta$ and $0<t<1$,
\begin{equation}\label{eq:app:chernoff}
 \chi(1)-(\chi-\eta)_+(1)
 \leq\htr(h_\chi^{1-t}h_\eta^t).
\end{equation}
Consequently, if $\rho,\sigma$ are states, $\lambda>0$, and $Q$ is
the support projection of $(\rho-\lambda\sigma)_+$, then
\begin{equation}\label{eq:app:test}
 \rho(1-Q)+\lambda\sigma(Q)
 \leq\lambda^t
 \langle\xi_\rho,\Delta_{\sigma,\rho}^{t}\xi_\rho\rangle.
\end{equation}
\end{lemma}
\begin{proof}
We use Ozawa's positive-power argument
\cite[Proposition~1.1]{JOPS} in Haagerup $L^p$ form.
Put $a=h_\chi$, $b=h_\eta$, and $c=b+(a-b)_+$. Then $c\geq a$ and
$c\geq b$. Positive-power monotonicity and positivity of the trace
pairing give, in succession,
\begin{align*}
 \htr a-\htr(a^{1-t}b^t)
 &=\htr\bigl(a^{1-t}(a^t-b^t)\bigr)\\
 &\leq\htr\bigl(a^{1-t}(c^t-b^t)\bigr)\\
 &\leq\htr\bigl(c^{1-t}(c^t-b^t)\bigr)\\
 &=\htr c-\htr(c^{1-t}b^t)\\
 &\leq\htr c-\htr b
 =\htr(a-b)_+.
\end{align*}
For clarity, the second inequality uses $c^t-b^t\geq0$ together with
$c^{1-t}-a^{1-t}\geq0$. The last uses
$c^{1-t}\geq b^{1-t}$, paired with $b^t$. All the products have finite
$L^1$ pairing by H\"older's inequality. This proves
\eqref{eq:app:chernoff}.

By the Jordan decomposition,
$(\rho-\lambda\sigma)(Q)=(\rho-\lambda\sigma)_+(1)$.
Hence
\[
 \rho(1-Q)+\lambda\sigma(Q)
 =1-(\rho-\lambda\sigma)_+(1).
\]
Apply the first part to $\chi=\rho$, $\eta=\lambda\sigma$, and then
use \eqref{eq:app:power} and
$h_{\lambda\sigma}^{t}=\lambda^t h_\sigma^t$.
This proves \eqref{eq:app:test}. For the general \vN\ algebra
inequality see Ogata \cite{Ogata2010}; its use for separation of states
is explained in \cite[Section~6]{JOPS}.
\end{proof}

We also need the Powers--St\o rmer estimate in standard form:
\begin{equation}\label{eq:app:ps}
 \|\xi_\chi-\xi_\eta\|^2\leq\|\chi-\eta\|.
\end{equation}
Here it follows from the same calculation, rather than being an
additional asymptotic input. At $t=1/2$, \eqref{eq:app:chernoff} gives
\[
 \langle\xi_\chi,\xi_\eta\rangle
 \geq\chi(1)-(\chi-\eta)_+(1).
\]
The inner product of these natural-cone vectors is real and nonnegative.
Using
$2(\chi-\eta)_+(1)=\|\chi-\eta\|+\chi(1)-\eta(1)$ and
expanding the squared vector norm proves \eqref{eq:app:ps}.

\subsection{A positive-part estimate gives an order bound}
\begin{lemma}[Order approximation]\label{lem:app:cutoff}
Let $\rho,\sigma$ be normal states, $\lambda>0$, and
$\delta=(\rho-\lambda\sigma)_+(1)<1$.
There is a normal state $\omega$ such that
\begin{equation}\label{eq:app:cutofforder}
 \omega\leq\frac{\lambda}{(1-\sqrt\delta)^2}\,\sigma
\end{equation}
and
\begin{equation}\label{eq:app:cutoffclose}
 P(\omega,\rho)\leq2\sqrt\delta,
 \qquad \|\omega-\rho\|\leq4\sqrt\delta.
\end{equation}
\end{lemma}
\begin{proof}
\step{Step 1: dominate both functionals by a common one.}
Set
\[
 A=\lambda\sigma,\qquad
 \gamma=(\rho-\lambda\sigma)_+,\qquad B=A+\gamma.
\]
Then $B\geq A$, $B\geq\rho$, and $\|B-A\|=\delta$.
The inequality $B\geq\rho$ permits an elementary commutant
construction: on $M\xi_B$ define
\[
 v'(x\xi_B)=x\xi_\rho.
\]
This is well defined and contractive, since
\[
 \|x\xi_\rho\|^2=\rho(x^*x)
 \leq B(x^*x)=\|x\xi_B\|^2.
\]
It therefore extends to $\overline{M\xi_B}$. Extend it by zero on
the orthogonal complement. Both $\overline{M\xi_B}$ and
$\overline{M\xi_\rho}$ reduce $M$, and the displayed rule intertwines
the action of $M$. Thus the extension is a contraction in $M'$.
In particular $v'\xi_B=\xi_\rho$.

\step{Step 2: apply the contraction to the smaller functional.}
Put $\eta=v'\xi_A$. For $x\geq0$, commutation with $x^{1/2}$ gives
\[
 \langle\eta,x\eta\rangle
 =\|v'x^{1/2}\xi_A\|^2
 \leq\|x^{1/2}\xi_A\|^2=A(x).
\]
Thus the functional represented by $\eta$ is bounded by
$A=\lambda\sigma$. Moreover \eqref{eq:app:ps} implies
\[
 \|\eta-\xi_\rho\|
 =\|v'(\xi_A-\xi_B)\|
 \leq\|\xi_A-\xi_B\|\leq\sqrt\delta.
\]

\step{Step 3: normalize and estimate the error.}
Let $r=\|\eta\|$. Since $\|\xi_\rho\|=1$, the reverse triangle
inequality gives $r\geq1-\sqrt\delta>0$. The unit vector
$\zeta=\eta/r$ represents a normal state $\omega$ satisfying
\eqref{eq:app:cutofforder}. Also
\[
 \|\zeta-\xi_\rho\|
 \leq\|\zeta-\eta\|+\|\eta-\xi_\rho\|
 =|1-r|+\|\eta-\xi_\rho\|\leq2\sqrt\delta.
\]
The definition of root fidelity gives
$F(\omega,\rho)\geq|\langle\zeta,\xi_\rho\rangle|$.
For unit vectors, $1-|\langle\zeta,\xi_\rho\rangle|^2
\leq\|\zeta-\xi_\rho\|^2$; hence the first bound in
\eqref{eq:app:cutoffclose} follows.
For the norm bound, expand the difference of the two vector states:
for $\|x\|\leq1$ its absolute value is at most
$2\|\zeta-\xi_\rho\|$. Taking the supremum proves the second bound.
The numerical constants are not important for the asymptotic exponents.
\end{proof}

\subsection{Why domination gives the required modular moments}
Write $e=\supp(\varphi)$, $\xi=\xi_\varphi$, and
$\Delta=\Delta_{\psi,\varphi}$. The relative Tomita operator has
$M\xi$ as a core on its initial support. On that core its formula is
$S_{\psi,\varphi}(x\xi)=e x^*\xi_\psi$. Therefore
\begin{equation}\label{eq:app:tomita}
 \|\Delta^{1/2}x\xi\|^2=\psi(x e x^*).
\end{equation}
The support $e$ is important when $\varphi$ is not faithful.

\begin{lemma}[Moment bounds]\label{lem:app:moments}
Let $\nu$ be the spectral probability measure of $\Delta$ in $\xi$.
Then $\nu(\{0\})=0$ and
\begin{equation}\label{eq:app:moments}
 \int_0^\infty r\,d\nu(r)\leq1,
 \qquad \int_0^\infty r^{-1}\,d\nu(r)\leq C.
\end{equation}
Consequently $s=S_M(\varphi\relsep\psi)$ is finite, and
\begin{equation}\label{eq:app:lambda}
 \Lambda(u):=\log\langle\xi,\Delta^{-u}\xi\rangle
 =\log\int_0^\infty e^{u(-\log r)}\,d\nu(r)
\end{equation}
is finite for $|u|\leq1$, with
\begin{equation}\label{eq:app:derivative}
 \Lambda(0)=0,\qquad \Lambda'(0)=s,
 \qquad 0\leq s\leq\log C.
\end{equation}
\end{lemma}
\begin{proof}
\step{Step 1: the positive and inverse moments.}
Taking $x=1$ in \eqref{eq:app:tomita} gives
$\int r\,d\nu(r)=\psi(e)\leq1$.
For the inverse moment, Cauchy--Schwarz for the state $\varphi$ and
$\varphi\leq C\psi$ yield
\begin{align*}
 |\langle\xi,x\xi\rangle|^2
 &=|\varphi(xe)|^2
 \leq\varphi(x e x^*)
 \leq C\psi(x e x^*)
 =C\|\Delta^{1/2}x\xi\|^2.
\end{align*}
Since $M\xi$ is a form core on the initial support, this extends to
$z$ in the form domain there:
\[
 |\langle\xi,z\rangle|^2\leq C\|\Delta^{1/2}z\|^2.
\]
Taking $z$ to be the projection of $\xi$ onto $\ker\Delta$ proves
$\nu(\{0\})=0$. More explicitly for the inverse moment, set
$z_k=\Delta^{-1}1_{[1/k,k]}(\Delta)\xi$. These vectors lie in the
form domain. If $I_k=\int_{[1/k,k]}r^{-1}\,d\nu(r)$, the last
inequality says $I_k^2\leq C I_k$, so $I_k\leq C$.
Monotone convergence proves the second part of
\eqref{eq:app:moments}.

\step{Step 2: differentiate the logarithmic moment.}
Put $L(r)=-\log r$. Since $e^{|L(r)|}\leq r+r^{-1}$,
\eqref{eq:app:moments} implies
$\int e^{|L|}\,d\nu\leq1+C$. Thus $L$ is integrable, and
differentiation under the integral at zero is justified by dominated
convergence (one may restrict to $|u|\leq1/2$).
It gives $\Lambda'(0)=\int L\,d\nu=s$ by
\eqref{eq:araki}. Jensen's inequality applied to $e^{-L}$ and $e^L$
gives respectively $e^{-s}\leq1$ and $e^s\leq C$.
This proves \eqref{eq:app:derivative}.
\end{proof}

Tensor products multiply relative modular operators and tensor their
representing vectors. Therefore the corresponding moment for the pair
$(\varphi_n,\psi_n)$ is
\begin{equation}\label{eq:app:tensor}
 \langle\xi_{\varphi_n},
 \Delta_{\psi_n,\varphi_n}^{-u}\xi_{\varphi_n}\rangle
 =e^{n\Lambda(u)}.
\end{equation}
Only the first-order expansion $\Lambda(u)=su+o(u)$ will be needed.

We record a useful consequence that controls the probability of an
arbitrary positive contraction. If $0\leq Q\leq1$, put
$p=\varphi(Q)$ and $q=\psi(Q)$. For $0<u<1$, the spectral
H\"older inequality gives
\begin{align*}
 \|\Delta^{u/2}Q\xi\|^2
 &\leq\|Q\xi\|^{2(1-u)}
       \|\Delta^{1/2}Q\xi\|^{2u}\\
 &\leq p^{1-u}q^u.
\end{align*}
Indeed $\|Q\xi\|^2=\varphi(Q^2)\leq p$, and
\eqref{eq:app:tomita} gives
$\|\Delta^{1/2}Q\xi\|^2=\psi(QeQ)\leq\psi(Q^2)\leq q$.
Since $\xi$ has no component in $\ker\Delta$, spectral calculus and
Cauchy--Schwarz now give
\[
 p^2
 =|\langle\Delta^{-u/2}\xi,\Delta^{u/2}Q\xi\rangle|^2
 \leq e^{\Lambda(u)}p^{1-u}q^u.
\]
For $p>0$ divide by $p^{1-u}$; for $p=0$ the resulting inequality
is automatic. We have proved
\begin{equation}\label{eq:app:binarymoment}
 p^{1+u}\leq e^{\Lambda(u)}q^u,
 \qquad 0<u<1.
\end{equation}
On the $n$th tensor power the exponent on the right is $n\Lambda(u)$.

\subsection{Asymptotic separation by positive contractions}
Recall that $\beta_{\varepsilon,n}$ is the infimum of $\psi_n(Q)$
over positive contractions with $\varphi_n(Q)\geq1-\varepsilon$.

\step{Direct bound: construct separating projections for each $R<s$.}
By \eqref{eq:app:derivative}, for every real $R<s$ there is a
sufficiently small $u\in(0,1)$ such that
$uR+\Lambda(-u)<0$. Apply \eqref{eq:app:test} with
$\lambda=e^{nR}$ and $t=u$, and let $Q_n$ be its positive-part
support projection. Then
\begin{equation}\label{eq:app:directtest}
 \varphi_n(1-Q_n)+e^{nR}\psi_n(Q_n)
 \leq e^{n(uR+\Lambda(-u))}\longrightarrow0.
\end{equation}
In particular $\varphi_n(Q_n)\to1$. For all large $n$ the
right side is at most $1$, so $\psi_n(Q_n)\leq e^{-nR}$.
For any fixed $\varepsilon>0$ these projections are eventually admissible;
hence
\[
 \liminf_{n\to\infty}-\frac1n\log\beta_{\varepsilon,n}\geq R.
\]
Letting $R\uparrow s$ proves the desired lower bound on the exponent.
Allowing negative $R$ makes this argument apply also when $s=0$.

\step{Converse bound: a moment estimate for positive contractions.}
For any admissible positive contraction, \eqref{eq:app:binarymoment} implies
\[
 \psi_n(Q)\geq
 (1-\varepsilon)^{(1+u)/u}e^{-n\Lambda(u)/u}.
\]
The lower bound is uniform in $Q$, so it also bounds the infimum
$\beta_{\varepsilon,n}$. This infimum is strictly positive; alternatively
$\varphi_n\leq C^n\psi_n$ gives
$\beta_{\varepsilon,n}\geq(1-\varepsilon)C^{-n}$ directly.
Take logarithms, divide by $n$, and first let $n\to\infty$:
\[
 \limsup_{n\to\infty}-\frac1n\log\beta_{\varepsilon,n}
 \leq\frac{\Lambda(u)}{u}.
\]
Now let $u\downarrow0$ and use \eqref{eq:app:derivative}.
Together with the direct bound this proves \eqref{eq:steinaep}.
The error parameter $\varepsilon\in(0,1)$ stays fixed throughout.

\subsection{Asymptotic logarithmic order bounds}
\label{app:upperorder}
\step{Upper bound: a positive-part estimate gives close, dominated states.}
Fix $R>s$. Choose $u\in(0,1)$ so small that
$\Lambda(u)-uR<0$, and put
\[
 \delta_n=(\varphi_n-e^{nR}\psi_n)_+(1).
\]
If $Q$ supports this positive part, write $p=\varphi_n(Q)$ and
$q=\psi_n(Q)$. When $\delta_n>0$ we have $p>e^{nR}q$.
Insert $q<e^{-nR}p$ into \eqref{eq:app:binarymoment} to obtain
$p\leq e^{n(\Lambda(u)-uR)}$. Since $\delta_n\leq p$, this proves
\begin{equation}\label{eq:app:hockey}
 0\leq\delta_n\leq e^{n(\Lambda(u)-uR)}\longrightarrow0.
\end{equation}
If $\delta_n=0$, the same bound holds without division by $p$.

Lemma~\ref{lem:app:cutoff} now gives, for all sufficiently large $n$,
\emph{normalized} states $\omega_n$ satisfying
\begin{align*}
 P(\omega_n,\varphi_n)&\leq2\sqrt{\delta_n}\longrightarrow0,\\
 \|\omega_n-\varphi_n\|&\leq4\sqrt{\delta_n}\longrightarrow0,\\
 \Dmax(\omega_n\relsep\psi_n)
 &\leq nR-2\log(1-\sqrt{\delta_n}).
\end{align*}
For fixed $\varepsilon>0$, these states are eventually feasible in
the definition of $\Dmax[,\varepsilon]_{M^{\otimes n}}$. The last correction divided by
$n$ tends to zero. Consequently
\[
 \limsup_{n\to\infty}\frac1n
 \Dmax[,\varepsilon]_{M^{\otimes n}}(\varphi_n\relsep\psi_n)\leq R.
\]
Let $R\downarrow s$. Notice the stronger conclusion along the way:
for \emph{each} $R>s$ we have normalized approximants converging in
norm and satisfying $\limsup_n n^{-1}\Dmax(\omega_n\relsep\psi_n)\leq R$. This is exactly
what is used in the diagonal choice of Lemma~\ref{lem:approximants}.

\step{Lower bound: approximation controls evaluation on $Q_n$.}
Fix $R<s$ and the projections $Q_n$ from
\eqref{eq:app:directtest}. Write
$a_n=\varphi_n(1-Q_n)\to0$; for all large $n$,
$\psi_n(Q_n)\leq e^{-nR}$.
Let $\alpha$ be any positive functional with $\alpha(1)\leq1$ and
$P(\alpha,\varphi_n)\leq\varepsilon$.
Because $\varphi_n(1)=1$, the extra term in $\overline F$ vanishes;
the condition on $P$ therefore gives
$F(\alpha,\varphi_n)\geq\sqrt{1-\varepsilon^2}$.

In any common vector representation with vectors $\eta$ for $\alpha$
and $\zeta$ for $\varphi_n$, decompose the inner product with $Q_n$
and $1-Q_n$. Cauchy--Schwarz gives
\begin{align*}
 |\langle\eta,\zeta\rangle|
 &\leq\sqrt{\alpha(Q_n)\varphi_n(Q_n)}
       +\sqrt{\alpha(1-Q_n)\varphi_n(1-Q_n)}\\
 &\leq\sqrt{\alpha(Q_n)}+\sqrt{a_n}.
\end{align*}
Taking the supremum over common representations yields the same bound
for $F$. Hence, for large $n$ and uniformly over all such $\alpha$,
\[
 \alpha(Q_n)\geq
 \bigl(\sqrt{1-\varepsilon^2}-\sqrt{a_n}\bigr)^2
 \geq c_\varepsilon:=\frac{1-\varepsilon^2}{4}>0.
\]
If also $\alpha\leq c\psi_n$, evaluation at $Q_n$ gives
$c\geq c_\varepsilon e^{nR}$. This holds for every feasible order
constant $c$ and every feasible $\alpha$. Taking both infima shows
\[
 \Dmax[,\varepsilon]_{M^{\otimes n}}(\varphi_n\relsep\psi_n)
 \geq nR+\log c_\varepsilon.
\]
Divide by $n$, take the lower limit, and let $R\uparrow s$.
This proves the lower bound and hence \eqref{eq:maxaep}.
Both assertions of Proposition~\ref{prop:aep} are now proved.


\end{document}